\documentclass{amsart}
\usepackage{amscd,amssymb, amsfonts}

\DeclareFontFamily{OMS}{rsfs}{\skewchar\font'60}
\DeclareFontShape{OMS}{rsfs}{m}{n}{<-5>rsfs5 <5-7>rsfs7 <7->rsfs10 }{}
\DeclareSymbolFont{rsfs}{OMS}{rsfs}{m}{n}
\DeclareSymbolFontAlphabet{\scr}{rsfs}

\newcommand{\sC}{\scr{C}}
\newcommand{\sD}{\scr{D}}

\newcommand{\sF}{\scr{F}}

\newcommand{\sP}{\scr{P}}

\renewcommand{\Im}{\operatorname{Im}}

\newcommand{\bb}{\begin{itemize}}

\newcommand{\ee}{\end{itemize}}
\newcommand{\ben}{\begin{enumerate}}
\newcommand{\een}{\end{enumerate}}

\newtheorem{thm}{Theorem}[section]
\newtheorem{theorem}{Theorem}[section]

\newtheorem{corollary}[thm]{Corollary}

\newtheorem{proposition}[thm]{Proposition}

\newtheorem{lemma}[thm]{Lemma}

\newtheorem*{thma}{Theorem 1.1}
\newtheorem*{thmb}{Theorem 1.2}

\theoremstyle{definition}

\newtheorem{definition}[thm]{Definition}

\theoremstyle{remark}

\newtheorem{remark}[thm]{Remark}

\newcommand{\Z}{\ensuremath{\mathbb{Z}}}

\newcommand{\ar}{\ensuremath{\operatorname{Ar}}}

\newcommand{\op}{\ensuremath{\operatorname{op}}}

\newcommand{\Top}{\ensuremath{\operatorname{Top}}}
\newcommand{\sd}{\ensuremath{\operatorname{Sd}}}
\newcommand{\id}{\ensuremath{\operatorname{id}}}

\input xy
\xyoption{all}
\CompileMatrices

\title{A discrete model of $S^1$-homotopy theory}
\author{Andrew J. Blumberg}
\address{Department of Mathematics \\ University of Chicago \\ 5734 S. University Ave. \\ Chicago, IL 60637} 
\email{blumberg@math.uchicago.edu}

\begin{document}

\makeatletter
\let\c@equation\c@thm
\makeatother
\numberwithin{equation}{section}

\begin{abstract}
We construct a discrete model of the homotopy theory of $S^1$-spaces.  We define a category $\sP$ with objects composed of a simplicial set and a cyclic set along with suitable compatibility data.  $\sP$ inherits a model structure from the model structures on the categories of simplicial sets and cyclic sets.  We then show that there is a Quillen equivalence between $\sP$ and the model category of $S^1$-spaces in which weak equivalences and fibrations are maps inducing weak equivalences and fibrations on passage to all fixed point sets.
\end{abstract}

\maketitle

\section{Introduction}

Simplicial techniques are often unavailable in the context of equivariant homotopy theory.  When $G$ is not a discrete group, simplicial $G$-sets do not provide a model for the homotopy theory of $G$-spaces.  The lack of an adequate replacement for simplicial sets is a substantial inconvenience.  Cyclic sets \cite{connes} provide a useful discrete model of a portion of $S^1$-homotopy theory.  Specifically, Spalinski \cite{spalinski:paper} (following Dwyer, Hopkins, and Kan \cite{dwyer-hopkins-kan}) constructs a model structure on cyclic sets which is Quillen equivalent to the model structure on $S^1$-spaces in which weak equivalences and fibrations are detected on passage to fixed point subspaces for finite groups.  However, since the $S^1$-fixed points of the geometric realization of a cyclic set must be discrete (\cite{fiedorowicz-gadja}, \cite{spalinski:paper}), it is unreasonable to expect a model structure on cyclic sets which will capture all of $S^1$-homotopy theory.

Restating this observation, the category of cyclic sets encodes all of $S^1$-homotopy theory except for the information detected by the $S^1$-fixed points.  A fundamental insight of Elmendorf \cite{elmendorf} is that the homotopy theory of $G$-spaces is equivalent to the homotopy theory of appropriate diagrams of fixed-point information.  See also Mandell and Scull \cite{mandell-scull} for a comprehensive modern discussion of this.  This suggests that a natural avenue of attack is to consider a category consisting of a cyclic set appropriately coupled (via compatibility data) with a simplicial set to represent the information at the $S^1$ fixed points.  Let $X$ be an $S^1$-space, and consider the following diagram:
\begin{equation}
\begin{CD}
X^{S^1} \times E\sF @>>> X \times E\sF \\
@VVV @. \\
X^{S^1}. @. \\
\end{CD}
\end{equation}
\\
Here $E\sF$ is the classifying space for the family of finite subgroups of $S^1$, the horizontal map is the inclusion and the vertical map is the projection.  The associated pushout is weakly equivalent to $X$.  This picture provides the inspiration for our construction.  We think of the cyclic set as akin to $X \times E\sF$, the simplicial set as $X^{S^1}$, and the compatibility data as the gluing along $X^{S^1} \times E\sF$. 

Given a simplicial set $A$ and a cyclic set $B$ we will describe the required compatibility in terms of a map $\nabla A \rightarrow B$, where $\nabla A$ is a homotopical cyclic approximation of $A$.  Specifically, we construct a functor $\nabla : {\bf S} \rightarrow {\bf S}^c$ which has the property that there is a natural map $|\nabla A|_c \rightarrow |A|_s$ which is a weak equivalence upon passage to all fixed point sets for finite subgroups of $S^1$.  The category $\sP$ of compatible pairs is an instance of a more general construction.

\newpage

\begin{definition}
Let $\sC$ and $\sD$ be categories and $F : \sC \rightarrow \sD$ a functor.  The category $\sC_F\sD$ has
\begin{enumerate}
\item{Objects specified by triples $(A, B, FA \rightarrow B)$ where $A$ is an object of $\sC$ and $B$ is an object of $\sD$.}
\item{Morphisms specified by maps $f_1 : A \rightarrow A^{\prime}$ and $f_2 : B \rightarrow B^{\prime}$ such that the following diagram commutes:
\begin{equation}
\begin{CD}
FA @>>> B \\
@V Ff_1 VV @VV f_2 V \\
FA^{\prime} @>>> B^{\prime}.\\  
\end{CD}
\end{equation}
}
\end{enumerate}
\end{definition}    

\begin{remark}
This is an example of a comma category \cite{maclane}.
\end{remark}

\begin{definition}
The category $\sP$ is the comma category ${\bf S}_{\nabla}{\bf S}^c$.
\end{definition}

When there are model structures on $\sC$ and $\sD$, there is an induced model structure on $\sC_F\sD$ for suitable functors $F$.

\begin{definition}
Let $\sC$ and $\sD$ be model categories.  A functor $F : \sC \rightarrow \sD$ is {\em Reedy admissible} if $F$ preserves colimits (e.g. $F$ is a left adjoint) and $F$ has the property that given a morphism $(A, B, FA \rightarrow B) \rightarrow (A^{\prime}, B^{\prime}, FA^{\prime} \rightarrow B^{\prime})$ in $\sC_F\sD$ such that $A \rightarrow A^{\prime}$ is a trivial cofibration in $\sC$ and $FA^{\prime} \cup_{FA} B \rightarrow B^{\prime}$ is a trivial cofibration in $\sD$ then $B \rightarrow B^{\prime}$ is a weak equivalence in $\sD$ (e.g. $F$ is a left Quillen functor).

\newdir{ (}{{}*!/-5pt/@^{(}}
\newdir{ ((}{{}*!/-15pt/@^{(}}

\begin{equation}
\xymatrix{
A \ar@{{ (}->}[dd]^\simeq & & FA \ar[rr]\ar[dd] & & B\ar[dd]\ar[dl] & & & & & B \ar[dd]^\simeq \\
& & & FA^{\prime} \cup_{FA} B \ar@{{ ((}->}[dr]^\simeq & & & \ar@2{->}[r] & &\\
A^{\prime} & & FA^{\prime} \ar[rr]\ar[ur] & & B^{\prime} & & & & & B^{\prime}\\\
}
\end{equation}
\end{definition}

\begin{theorem}
Let $\sC$ and $\sD$ be model categories and $F : \sC \rightarrow \sD$ be a Reedy admissible functor.  Then $\sC_F\sD$ admits a model structure.  A map $(A,B,FA \rightarrow B) \rightarrow (A^{\prime},B^{\prime},FA^{\prime} \rightarrow B^{\prime})$ is  

\begin{enumerate}
\item{a weak equivalence if $A \rightarrow A^{\prime}$ is a weak equivalence in $\sC$ and $B \rightarrow B^{\prime}$ is a weak equivalence in $\sD$,} 
\item{a fibration if $A \rightarrow A^{\prime}$ is a fibration in $\sC$ and $B \rightarrow B^{\prime}$ is a fibration in $\sD$,}  
\item{a cofibration if $A \rightarrow A^{\prime}$ is a cofibration in $\sC$ and $FA^{\prime} \cup_{FA} B \rightarrow B^{\prime}$ is a cofibration in $\sD$.}
\end{enumerate}
\end{theorem}

We use this theorem to obtain the model structure on $\sP$.

\begin{lemma}
The functor $\nabla : {\bf S} \rightarrow {\bf S}^c$ is Reedy admissible.
\end{lemma}

\begin{corollary}
There is a model structure on $\sP$ inherited from the model structures on ${\bf S}$ and ${\bf S}^c$.
\end{corollary}

There is an adjunction specified by functors $L : \sP \rightarrow \Top^{S^1}$ and $R : \Top^{S^1} \rightarrow \sP$.  We will provide more details about these functors later.

\begin{definition}
The functor $L: \sP \rightarrow \Top^{S^1}$ takes a triple $(A, B, \nabla A \rightarrow B)$, to the pushout in the diagram:

\begin{equation}
\begin{CD}
|\nabla A|_c @>>> |B|_c \\
@VVV @VVV \\
|A|_s @>>> X. \\
\end{CD}
\end{equation}

The functor $R:\Top^{S^1} \rightarrow \sP$ takes $X$ to the triple $(S(X^{S^1}), S_c(X), \nabla S(X^{S^1}) \rightarrow S_c(X))$.  The map $\xi : \nabla S(X^{S^1}) \rightarrow S_c(X)$ is the adjoint of the composite 
\begin{equation}|\nabla S(X^{S^1})|_c \rightarrow |S(X^{S^1})|_s \rightarrow X^{S^1} \hookrightarrow X.
\end{equation}
\end{definition}

Recall that there is a model structure on $\Top^{S^1}$ given by defining a map $f:X \rightarrow Y$ to be a weak equivalence if $f^H : X^H \rightarrow Y^H$ is a weak equivalence for all $H \subset S^1$, a fibration if $f^H : X^H \rightarrow Y^H$ is a fibration for all $H \subset S^1$, and a cofibration if it has the left-lifting property with respect to acyclic fibrations \cite{spalinski:paper}. 

Here is the main theorem of the paper.

\begin{theorem}
The functors $L$ and $R$ specify a Quillen equivalence between $\sP$ with the model structure given by Theorem 1.1 and $\Top^{S^1}$ with the model structure described above.
\end{theorem}

The problem of obtaining a discrete model for $S^1$-spaces was raised by Voevodsky in a 2002 e-mail to May \cite{may-voevodsky}.

The rest of the paper is organized as follows: 
\begin{enumerate}
\item{A brief review of simplicial and cyclic sets.}
\item{A review of the model structures on $S^1$ spaces.}
\item{Definition of $\nabla$ and demonstration that it is Reedy admissible.}
\item{The adjunction between $\sP$ and $\Top^{S^1}$.}
\item{Proof of Theorem 1.2.}
\item{Proof of Theorem 1.1.}
\item{Appendix: calculations from Spalinski's thesis \cite{spalinski:thesis}.}
\end{enumerate}

\section{A review of cyclic sets}

We give a very succinct review of cyclic sets.  Good references for readers unfamiliar with the category are \cite{spalinski:paper, dwyer-hopkins-kan, fiedorowicz-gadja, connes}.  A cyclic set can be regarded as a simplicial set with extra data, namely an action of $\Z/(n+1)$ on the n-simplices which is compatible with the face and degeneracy operators.  Define the cyclic category $\Lambda^{\op}$ to have the same objects as the category $\Delta^{\op}$ and the same generating morphisms along with an extra degeneracy $s_{n+1} : [n] \rightarrow [n+1]$ and the ``cyclic relations'' $(d_0 s_{n+1})^{n+1} = \id$.  Define $t_n = d_0 s_{n+1} : [n] \rightarrow [n]$.  Every morphism in $\Lambda^{\op}$ can be written as a composite $\phi = STD$ of a composite $S$ of degeneracy operators, a power $T$ of $t_n$ for some $n$, and a composite $D$ of boundary operators \cite{spalinski:paper}.  For further discussion of the properties of $\Lambda^{\op}$ (e.g. $\Lambda^{\op}$ is self dual) see \cite{connes} or \cite{elmendorf2}. 

Cyclic sets are contravariant functors from the category $\Lambda^{\op}$ to sets.  The category of cyclic sets will be denoted ${\bf S}^c$.  As in the theory of simplicial sets, the represented cyclic sets $\Lambda[n] = \hom_{\Lambda}(-,n)$ play an important role.  The geometric realization of the underlying simplicial set of a cyclic set admits a natural $S^1$-action.  The geometric realization, regarded as a functor from cyclic sets to $S^1$-spaces, will be denoted $|-|_c$.  In particular, $|\Lambda[n]|_c \cong S^1 \times |\Delta[n]|$, with $S^1$ acting on the product by rotation on the first coordinate, where $\Delta[n]$ is the represented simplicial set with trivial action.  By manipulation of coends one obtains $|X|_c = X \otimes_{\Lambda^{\op}} |\Lambda|_c$.  The adjoint to the realization is the ``cyclic singular functor'' $S_c$ defined to have $n$-simplices $\hom_{S^1}(|\Lambda[n]|_c,X)$ \cite{dwyer-hopkins-kan}.  Here the cyclic structure is obtained by regarding $|\Lambda[n]|$ as homeomorphic to $S^1 \times |\Delta[n]|$, where the action of $t_n$ permutes the coordinates of a point in $|\Delta[n]|$ and rotates $S^1$ by $e^{2\pi i / (n+1)}$.

Now consider the subgroup $\Z / (r) \subset S^1$.  Given a cyclic set, we can apply the subdivision functor $\sd_r$ to the underlying simplicial set \cite{bokstedt-hsiang-madsen}.  This has a natural simplicial action of $\Z / (r)$ induced from the cyclic structure, and so we can define a composite functor $\Phi_r$ which takes a cyclic set $X$ to the simplicial set $(\sd_r X)^{\Z / (r)}$.  There is a homeomorphism $|\Phi_r(X)|_s \cong |X|^{\Z / (r)}$.  By Freyd's adjoint functor theorem, $\Phi_r$ has an adjoint $\Psi_r$.  It is useful to describe $\Psi_r$ more concretely and so we reproduce calculations of Spalinski \cite{spalinski:thesis} in the appendix.

The functors $\Phi_r$ are used to prove the following result.

\begin{lemma}[Spalinski \cite{spalinski:paper}]\label{counit}
The counit of the adjunction between $|-|_c$ and $S_c(-)$ induces weak equivalences on passage to fixed point spaces for finite subgroups of $S^1$.
\end{lemma}

\section{A review of model structures on simplicial sets, cyclic sets, and $\Top^{S^1}$}

We briefly review the model structures on ${\bf S}$, ${\bf S}^c$, and $\Top^{S^1}$.

\begin{theorem} 
There is a model structure on simplicial sets in which a map is
\begin{enumerate}
\item{a fibration if it is a Kan fibration,}
\item{a weak equivalence if the induced map on passage to geometric realization is an equivalence,}
\item{a cofibration if it is an injection.}
\end{enumerate}
\end{theorem}

\begin{theorem} 
For any family $\sF$ of subgroups of $S^1$, there is a model structure on $\Top^{S^1}$ in which a map $f: X \rightarrow Y$ is
\begin{enumerate}
\item{a fibration if the induced maps $f^H : X^H \rightarrow Y^H$ are fibrations for all $H \in \sF$,}
\item{a weak equivalence if the induced maps $f^H : X^H \rightarrow Y^H$ are weak equivalences for all $H \in \sF$,}
\item{a cofibration if it has the left-lifting property with respect to trivial fibrations.}
\end{enumerate}
In particular, this holds when $\sF$ is the family of all subgroups of $S^1$ and when $\sF$ is the family of all finite subgroups of $S^1$.
\end{theorem}

\begin{theorem}[Spalinski \cite{spalinski:paper}] 
There is a model structure on cyclic sets in which a map is  
\begin{enumerate}
\item{a fibration if $\Phi_r(f)$ is a fibration of simplicial sets for all $r \geq 1$,}
\item{a weak equivalence if $\Phi_r(f)$ is a weak equivalence of simplicial sets for all $r \geq 1$,}
\item{a cofibration if it has the left-lifting property with respect to trivial fibrations.}
\end{enumerate}
\end{theorem}

\begin{remark}
Spalinski \cite{spalinski:paper} also shows that cofibrations can be characterized as retracts of transfinite composites of pushouts of coproducts of maps $\Psi_r(\partial \Delta[k]) \rightarrow \Psi_r(\Delta[k])$.
\end{remark}

The homotopy theory of cyclic sets is the same as the homotopy theory of $\Top^{S^1}$ with respect to the family of finite subgroups of $S^1$.

\begin{theorem}[Spalinski \cite{spalinski:paper}]
The cyclic realization functor and the cyclic singular functor induce a Quillen equivalence between $\Top^{S^1}$ with the model structure in which $\sF$ is the family of finite subgroups of $S^1$ and the category of cyclic sets with the model structure described above.
\end{theorem}

\section{The functor $\nabla$}

We shall construct a functor $\nabla : {\bf S} \rightarrow {\bf S}^c$ such that $|X|_{s}^H \simeq |\nabla X|_{c}^H$ for all finite $H \subset S^1$.  One's first guess is that $\nabla$ ought to be the left adjoint to the forgetful functor which assigns to a cyclic set its underlying simplicial set (Kan extension).  However, this is the free cyclic set associated to the underlying simplicial set \cite{burghelea-fiedorowicz}, and does not have the properties we need.

Another obvious guess is to define $\nabla X = S_{c}(|X|_{s})$.  By Lemma \ref{counit}, we know the counit provides a map $|S_c(|X|_{s})|_c \rightarrow |X|_s$ which is an equivalence on passage to all finite subgroups.  Unfortunately, as a composite of a left adjoint and a right adjoint, this functor has rather unpleasant properties.  For instance, it preserves neither colimits nor limits.

We want a functor from simplicial sets to cyclic sets which is a left adjoint and so preserves colimits.  All such functors arise from cosimplicial cyclic sets.  In fact, there is an equivalence between the category of cosimplicial objects in $\sC$ and adjunctions from simplicial sets to $\sC$ for categories $\sC$ with all small colimits \cite[3.1.5]{hovey}.

\begin{definition}
Set $\nabla_n = S_c(|\Delta[n]|)$.  Then $\nabla_*$ is a cosimplicial cyclic set and so we can define a functor $\nabla : {\bf S} \rightarrow {\bf S}^c$ by letting $\nabla X = X \otimes_{\Delta^{\op}} \nabla_*$.
The functor $\nabla$ has the right adjoint $A : {\bf S}^c \rightarrow {\bf S}$ specified by
$A(Y)_n = \hom_{{\bf S}^c}(\nabla_n, Y)$.
\end{definition}  

We will repeatedly use the following result, which we quote from \cite{mandell-may}.

\begin{lemma}\label{fixedpushout}
The functor $(-)^G$ on based $G$-spaces preserves pushouts of diagrams one leg of which is a closed inclusion.
\end{lemma}

\begin{lemma}
There is a natural map $\zeta : |\nabla X|_c \rightarrow |X|_s$ which induces weak equivalences on passage to fixed point subspaces for all finite subgroups of $S^1$.
\end{lemma}

\begin{proof}
By construction, the counit map $\gamma_n : |\nabla_n|_c \rightarrow |\Delta[n]|_s$ induces weak equivalences on passage to all fixed point subspaces for finite subgroups of $S^1$.  Define $\zeta$ to be the following map:
\begin{eqnarray*}
|\nabla X|_c = ((X \otimes_{\Delta^{\op}} \nabla) \otimes_{\Delta^{\op}} |\Delta|_s) = (X \otimes_{\Delta^{\op}} (\nabla \otimes_{\Delta^{\op}} |\Delta|_s)) = X \otimes_{\Delta^{\op}} |\nabla|_c \longrightarrow X \otimes_{\Delta^{\op}} |\Delta|_s = |X|_s.
\end{eqnarray*}
Both the domain and the codomain can be regarded as a succession of pushouts with one leg a cofibration.  Therefore the fixed-point functor commutes with each of these coends by Lemma \ref{fixedpushout} and so $\zeta$ induces weak equivalences on passage to fixed subspaces.  
\end{proof}

\begin{remark}
The essential aspect of the $\nabla_n$ is that they come equipped with maps from $|\nabla_n|_c$ to $|\Delta[n]|_s$ which induce weak equivalences on passage to fixed point subspaces for all finite subgroups of $S^1$.  Any other cosimplicial cyclic set which had this property would suffice for our purposes.  One might prefer a functorial cofibrant approximation of $\nabla_*$.  Alternatively, as the singular construction we give is rather bloated, we expect that other explicit models of $\nabla_n$ may well be preferable for specific applications.  
\end{remark}

To use Theorem 1.1 to show that there is a model structure on $\sP$, we must verify that $\nabla$ is Reedy admissible.  By construction, $\nabla$ is a left adjoint and so preserves colimits.

\begin{lemma}\label{nablaadmit}
Given a map $(A, B, \nabla A \rightarrow B) \rightarrow (A^{\prime}, B^{\prime}, \nabla A^{\prime} \rightarrow B^{\prime})$ in $\sP$ such that $A \rightarrow A^{\prime}$ is a trivial cofibration and $\nabla A^{\prime} \cup_{\nabla A} B \rightarrow B^{\prime}$ is a trivial cofibration, the map $B \rightarrow B^{\prime}$ is a weak equivalence.  Therefore $\nabla$ is Reedy admissible.
\end{lemma}

\begin{proof}
Since $B \rightarrow B^{\prime}$ is the composite
\begin{equation}B \rightarrow \nabla A^{\prime} \cup_{\nabla A} B \rightarrow B^{\prime}\end{equation}and $\nabla A^{\prime} \cup_{\nabla A} B \rightarrow B^{\prime}$ is a weak equivalence by hypothesis, it suffices to show that $B \rightarrow \nabla A^{\prime} \cup_{\nabla A} B$ is a weak equivalence.  This is equivalent to showing that $|B|_c^H \rightarrow |\nabla A^{\prime} \cup_{\nabla A} B|_c^H$ is a weak equivalence of spaces for all finite $H \subset S^1$.  Since geometric realization is a colimit, $|\nabla A^{\prime} \cup_{\nabla A} B|_c^H$ is isomorphic to $(|\nabla A^{\prime}|_c \cup_{|\nabla A|_c} |B|_c)^H$.  Since $A \rightarrow A^{\prime}$ is a cofibration of simplicial sets (and hence an inclusion), $|\nabla A|_c \rightarrow |\nabla A^{\prime}|_c$ is a closed inclusion.  Therefore by Lemma \ref{fixedpushout} the fixed-point functor commutes with the pushout, and so $(|\nabla A^{\prime}|_c \cup_{|\nabla A|_c} |B|_c)^H$ is equivalent to $|\nabla A^{\prime}|_c^H \cup_{|\nabla A|_c^H} |B|_c^H$.  Since $|\nabla A|_c^H \rightarrow |\nabla A^{\prime}|_c^H$ is a trivial cofibration when $A \rightarrow A^{\prime}$ is a trivial cofibration, $|B|^H_c \rightarrow |\nabla A^{\prime}|_c^H \cup_{|\nabla A|_c^H} |B|_c^H$ is the pushout of a trivial cofibration and so is itself a trivial cofibration.       
\end{proof}

\begin{corollary}\label{pmodel}
There is a model structure on $\sP$ in which a map is
\begin{enumerate}
\item{a weak equivalence if $A \rightarrow A^{\prime}$ is a weak equivalence of simplicial sets and $B \rightarrow B^{\prime}$ is a weak equivalence of cyclic sets,}
\item{a fibration if $A \rightarrow A^{\prime}$ is a fibration of simplicial sets and $B \rightarrow B^{\prime}$ is a fibration of cyclic sets,}
\item{a cofibration if $A \rightarrow A^{\prime}$ is a cofibration of simplicial sets and $\nabla A^{\prime} \cup_{\nabla A} B \rightarrow B^{\prime}$ is a cofibration of cyclic sets.}
\end{enumerate}
\end{corollary}

\begin{proof}
This follows immediately from Lemma \ref{nablaadmit} and Theorem 1.1.
\end{proof}

\section{The adjunction between $\sP$ and $\Top^{S^1}$}

There are natural functors from $\sP$ to $\Top^{S^1}$ and from $\Top^{S^1}$ to $\sP$ defined as follows.

\begin{lemma}\label{mapcommute}
Given a morphism in $\sP$ from $(A, B, \nabla A \rightarrow B)$ to $(A^{\prime}, B^{\prime}, \nabla A^{\prime} \rightarrow B^{\prime})$, the induced diagram
\begin{equation}
\begin{CD}
|\nabla A|_c @>\zeta >> |A|_s \\
@VVV @VVV \\
|\nabla A^{\prime}|_c @>\zeta >> |A^{\prime}|_s \\
\end{CD}
\end{equation}
is commutative.
\end{lemma}

\begin{proof}
Rewriting the diagram as follows
\begin{equation}
\begin{CD}
A \otimes_{\Delta^{\op}} |\nabla| @>>> A \otimes_{\Delta^{\op}} |\Delta| \\
@VVV @VVV \\
A^{\prime} \otimes_{\Delta^{\op}} |\nabla| @>>> A^{\prime} \otimes_{\Delta^{\op}} |\Delta| \\
\end{CD}
\end{equation}
makes the commutativity apparent.
\end{proof} 

\begin{definition}
The functor $L: \sP \rightarrow \Top^{S^1}$ takes a triple $(A, B, \nabla A \rightarrow B)$, to the pushout in the diagram:

\begin{equation}
\begin{CD}
|\nabla A|_c @>>> |B|_c \\
@V\zeta VV @VVV \\
|A|_s @>>> X. \\
\end{CD}
\end{equation}
\\
A morphism $(A,B,\nabla A \rightarrow B) \rightarrow (A^{\prime}, B^{\prime}, \nabla A^{\prime} \rightarrow B^{\prime})$ induces a commutative diagram:
\begin{equation}
\begin{CD}
|A|_s @<<< |\nabla A|_c @>>> |B|_c \\
@VVV @VVV @VVV \\
|A^{\prime}|_s @<<< |\nabla A^{\prime}|_c @>>> |B^{\prime}|_c.\\
\end{CD}
\end{equation}
The lefthand square commutes by the preceding lemma and the righthand square commutes because of the definition of a morphism.  Therefore there is an induced map of pushouts, which specifies the action of $L$ on morphisms. 
\end{definition}

\begin{lemma}
A morphism $X \rightarrow Y$ in $\Top^{S^1}$ induces a commutative diagram:
\begin{equation}
\begin{CD}
\nabla S(X^{S^1}) @>\xi >> S_c(X) \\
@VVV @VVV \\
\nabla S(Y^{S^1}) @>\xi >> S_c(Y).\\
\end{CD}
\end{equation}
\end{lemma}

\begin{proof}
This diagram commutes if and only if the adjoint diagram
\begin{equation}
\begin{CD}
|\nabla S(X^{S^1})|_c @>>> X \\
@VVV @VVV \\
|\nabla S(Y^{S^1})|_c @>>> Y\\
\end{CD}
\end{equation}
commutes.  The latter diagram can be written as the composite: 
\begin{equation}
\begin{CD} 
|\nabla S(X^{S^1})|_c @>>> |S(X^{S^1})|_s @>>> X^{S^1} @>>> X \\
@VVV @VVV @VVV @VVV \\
|\nabla S(Y^{S^1})|_c @>>> |S(Y^{S^1})|_s @>>> Y^{S^1} @>>> Y.\\
\end{CD}
\end{equation}
Here the lefthand square commutes by Lemma \ref{mapcommute}.  The middle square commutes by the naturality of the counit.  The righthand square commutes trivially.  Therefore the original diagram commutes. 
\end{proof}

\begin{definition}
The functor $R:\Top^{S^1} \rightarrow \sP$ takes $X$ to the triple $(S(X^{S^1}), S_c(X), \nabla S(X^{S^1}) \rightarrow S_c(X))$.  The map $\nabla S(X^{S^1}) \rightarrow S_c(X)$ is the adjoint of the composite:
\begin{equation}|\nabla S(X^{S^1})|_c \rightarrow |S(X^{S^1})|_s \rightarrow X^{S^1} \hookrightarrow X\end{equation}
A map $X \rightarrow Y$ in $\Top^{S^1}$ induces maps $S(X^{S^1}) \rightarrow S(Y^{S^1})$ and $S_c(X) \rightarrow S_c(Y)$ by functoriality.  By the preceding lemma, these maps fit into a commutative diagram:
\begin{equation}
\begin{CD}
\nabla S(X^{S^1}) @>>> S_c(X) \\
@VVV @VVV \\
\nabla S(Y^{S^1}) @>>> S_c(Y).\\
\end{CD}
\end{equation}  
\end{definition}

We think of $L$ as a realization functor and $R$ as a singular functor.

\begin{proposition}
The functors
\begin{equation}L : \sP \rightleftarrows \Top^{S^1} : R\end{equation}
form an adjoint pair.
\end{proposition}

\begin{proof}
Given a map $|A|_s \cup_{|\nabla A|_c} |B|_c \rightarrow X$, we must show that there is a unique corresponding map $(A, B, \nabla A \rightarrow B) \rightarrow (S(X^{S^1}), S_c(X), \nabla S(X^{S^1}) \rightarrow S_c(X))$.  We clearly get unique maps $A \rightarrow S(X)$ and $B \rightarrow S_c(X)$ as adjoints to the maps $|A|_s \rightarrow X$ and $|B|_c \rightarrow X$ induced by the map from the pushout.  It suffices to verify that the compatibility imposed by the pushout square is equivalent to the compatibility condition for a morphism in $\sP$.

So consider the square induced by our adjoint maps:
\begin{equation}
\tag{*}
\begin{CD}
\nabla A @>>> B\\
@VVV @VVV \\
\nabla S(X^{S^1}) @>>> S_c(X). \\
\end{CD}
\end{equation}
We must show that it commutes.  Now, the map $|A|_s \cup_{|\nabla A|_c} |B|_c \rightarrow X$ provides us with a commuting square: 
\begin{equation}
\begin{CD}
|\nabla A|_c @>>> |B|_c \\
@VVV @VVV \\
|A|_s @>>> X. \\
\end{CD}
\end{equation}
Such squares are in bijective correspondence with commuting squares:
\begin{equation}
\tag{**}
\begin{CD}
\nabla A @>>> B \\
@VVV @VVV \\
S_c(|A|_s) @>>> S_c(X).\\
\end{CD}
\end{equation}
The two composites $\nabla A \rightarrow B \rightarrow S_c(X)$ are the same.  Therefore to verify the correspondence of the compatibility conditions it suffices to show that the maps \begin{equation}(*) \quad \nabla A \rightarrow \nabla S(X^{S^1}) \rightarrow S_c(X) \qquad\qquad (**) \quad \nabla A \rightarrow S_c(|A|_s) \rightarrow S_c(X)\end{equation} are identical.  We do this by explicitly chasing elements around these two paths.  Start with the map $g : |A|_s \rightarrow X$.  Regarding $|A|_s$ as the coend $A \otimes_{\Delta^{\op}} |\Delta|$, we view $g$ as taking $(a, \delta)$ to $g(a,\delta)$ and its adjoint as taking $a$ to the map $(\delta \rightarrow g(a,\delta))$.

So let's unwind the two maps.  The map 
\begin{equation}\tag{**} \nabla A \rightarrow S_c(|A|_s) \rightarrow S_c(X)\end{equation}
is the composite 
\begin{equation}\nabla A \rightarrow S_c(|A|_s) \rightarrow S_c(|A|_s) \rightarrow S_c(X).\end{equation}
The first constituent map is adjoint to the map $|\nabla A|_c \rightarrow |A|_s$ which we defined as 
\begin{equation}A \otimes_{\Delta^{\op}} |\nabla|_c \rightarrow A \otimes_{\Delta^{\op}} |\Delta|_s\end{equation}
via the map $\gamma : |\nabla|_c \rightarrow |\Delta|_s$.  In order to calculate the adjoint map, we write the first coend as 
\begin{equation}(A \otimes_{\Delta^{\op}} \nabla) \otimes_{\Lambda^{\op}} |\Lambda|_c \rightarrow A \otimes_{\Delta^{\op}} |\Delta|_s\end{equation}
where the map takes $((a, \nu), \lambda)$ to $(a, \gamma(\nu,\lambda))$.  Then the adjoint is the map 
\begin{equation}\lambda \rightarrow ((a,\nu) \rightarrow (a,\gamma(\nu,\lambda)).\end{equation}
Next, we have the map $S_c(|A|_s) \rightarrow S_c(X)$ which is obtained by applying $S_c$ to the map $g : |A|_s \rightarrow X$.  That is, the induced map takes the map $\lambda \rightarrow (a,\delta)$ to the map $\lambda \rightarrow g(a,\delta)$.  Finally, the composite is 
\begin{equation}\tag{**} (a,\nu) \rightarrow (\lambda \rightarrow g(a,\gamma(\nu,\lambda)))\end{equation}

On the other hand, we can decompose the map 
\begin{equation}\tag{*} \nabla A \rightarrow \nabla S(X^{S^1}) \rightarrow S_c(X)\end{equation} as the composite 
\begin{equation}\nabla A \rightarrow \nabla S(X^{S^1}) \rightarrow \nabla S(X^{S^1}) \rightarrow S_c(X).\end{equation}
The first constituent map is obtained by applying $\nabla$ to the map $A \rightarrow S(X^{S^1})$ adjoint to $g$.  Explicitly, this is 
\begin{equation}(a,\nu) \rightarrow ((a \rightarrow (\delta \rightarrow g(a,\delta))),\nu).\end{equation}
The second map is the adjoint to the map $|\nabla S(X^{S^1})|_c \rightarrow X$, which decomposes as the composite 
\begin{equation}|\nabla S(X^{S^1})|_c \rightarrow |S(X^{S^1})|_s \rightarrow X^{S^1} \rightarrow X\end{equation}
that takes $(h,\nu,\lambda)$ to $(h,\gamma(\nu,\lambda))$ and then to $h(\gamma(\nu,\lambda))$.  The adjoint can be written as: 
\begin{equation}(h, \nu) \rightarrow (\lambda \rightarrow h(\gamma(\nu,\lambda)).\end{equation}
Composing, we have 
\begin{equation}\tag{*} (a, \nu) \rightarrow (\lambda \rightarrow g(a,\gamma(\nu,\lambda))).\end{equation}
\end{proof}

\section{Proof of Theorem 1.2}

The functors $L$ and $R$ are compatible with our model structures.

\begin{lemma}
Let $\sP$ have the model structure described in Corollary \ref{pmodel} and $\Top^{S^1}$ have the model structure generated by the family of all subgroups of $S^1$.  Then the adjoint functors $L$ and $R$ form a Quillen adjunction.
\end{lemma}

\begin{proof} 
It suffices to show that $R$ preserves fibrations and trivial fibrations.  If $X \rightarrow Y$ is a fibration or a trivial fibration, then $S(X) \rightarrow S(Y)$ and $S_c(X) \rightarrow S_c(Y)$ are as well since $S_c(-)$ and $S(-)$ are themselves right Quillen adjoints.
\end{proof}

\begin{remark}
In fact, $R$ preserves weak equivalences since both $S(-)$ and $S_c(-)$ preserve weak equivalences.
\end{remark}

Now, one potential problem with this model for $\Top^{S^1}$ is that while cyclic sets don't capture ``useful'' data at the $S^1$ fixed points, they do have some information there which might corrupt the data encoded in the simplicial set.  In fact, it isn't in general the case that the counit $|S(X^{S^1})|_s \cup_{|\nabla S(X^{S^1})|_c} |S_c(X)|_c \rightarrow X$ is a weak equivalence of $S^1$-spaces.  However, the following lemmas show that this map is an equivalence once we pass to cofibrant approximations.  Observe that $(A,B, \nabla A \rightarrow B)$ cofibrant implies that $A$ is cofibrant and $\nabla A \rightarrow B$ is a cofibration.

\begin{lemma}\label{fixfreecof}
If $X \rightarrow Y$ is cofibration of cyclic sets, then the induced map $|X|_c^{S^1} \rightarrow |Y|_c^{S^1}$ is a homeomorphism.
\end{lemma}

\begin{proof}
As noted previously, a cofibration of cyclic sets is a retract of a relative cell complex with respect to the family $\Psi_r(\partial \Delta[k]) \rightarrow \Psi_r(\Delta[k])$.  A retract of a homeomorphism is a homeomorphism.  Thus, it will suffice to observe that the domains and codomains of these generating cofibrations have no $S^1$-fixed points, as the fixed-point functor commutes with these pushouts after passage to cyclic realization by Lemma \ref{fixedpushout}.  But this is true by an explicit calculation \cite{spalinski:thesis} which we reproduce in the appendix.
\end{proof}

\begin{corollary}
If $X$ is a cofibrant cyclic set, then $|X|_c^{S^1}$ is empty.
\end{corollary}

Lemma \ref{fixfreecof} enables us to show that for cofibrant objects in $\sP$ the gluing behaves properly.

\begin{lemma}\label{keylemma}
Let $(A, B, \nabla A \rightarrow B)$ be a cofibrant object in $\sP$ and define $Z = |A|_s \cup_{|\nabla A|_c} |B|_c$.  Then $Z^{S^1} \simeq |A|_s$ and for all finite $H \subset S^1$, $Z^H \simeq |B|_c^H$.
\end{lemma}

\begin{proof}
Observe that by Lemma \ref{fixedpushout}, passage to fixed points commutes with the pushout since $\nabla A \rightarrow B$ is a cofibration.  First consider the $S^1$-fixed points.  We have a map $|A|_s \rightarrow Z^{S^1}$ induced by the pushout.  Since $\nabla A \rightarrow B$ is a cofibration, Lemma \ref{fixfreecof} tells us that $|\nabla A|_c^{S^1} \cong |B|_c^{S^1}$.  But this immediately implies that $(|A|_s \cup_{|\nabla A|_c} |B|_c)^{S^1} \cong |A|_s$.  Now consider a finite subgroup $H \subset S^1$.  Since $|\nabla A|_c^H \simeq |A|_s^H$ and $|\nabla A|_c^H \rightarrow |B|_c^H$ is a cofibration, $|B|_c^H \rightarrow Z^H$ is the pushout along a cofibration of a weak equivalence.  Therefore $|B|_c^H \rightarrow Z^H$ is a weak equivalence since $\Top$ is proper.
\end{proof}

\begin{thmb}
The functors $L$ and $R$ specify a Quillen equivalence between $\sP$ with the model structure given by Theorem 1.1 and $\Top^{S^1}$ with the model structure in which $\sF$ is the family of all subgroups of $S^1$. 
\end{thmb}

\begin{proof}
We must show that given a cofibrant object $(A, B, \nabla A \rightarrow B)$ in $\sP$ and a fibrant $S^1$-space $X$, a map $(A,B,\nabla A \rightarrow B) \rightarrow RX$ is a weak equivalence if and only if the adjoint $L(A,B,\nabla A \rightarrow B) \rightarrow X$ is a weak equivalence.  Writing out the functors, we need to show that \begin{equation}(A,B,\nabla A \rightarrow B) \rightarrow (S(X^{S^1}), S_c(X), \nabla S(X^{S^1}) \rightarrow S_c(X))\end{equation} is a weak equivalence if and only if \begin{equation}|A|_s \cup_{|\nabla A|_c} |B|_c \rightarrow X\end{equation} is a weak equivalence.

So assume that $|A|_s \cup_{|\nabla A|_c} |B|_c \rightarrow X$ is a weak equivalence.  This implies that the induced map 
\begin{equation}(|A|_s \cup_{|\nabla A|_c} |B|_c)^{S^1} \rightarrow X^{S^1}\end{equation}
is a weak equivalence.  Furthermore, by Lemma \ref{keylemma} the map 
\begin{equation}|A|_s \rightarrow (|A|_s \cup_{|\nabla A|_c} |B|_c)^{S^1}\end{equation}
is a weak equivalence, and so the composition is a weak equivalence.  This implies that the adjoint $A \rightarrow S(X^{S^1})$ is a weak equivalence of simplicial sets.  Similarly, for any finite $H \subset S^1$ the assumption implies that the induced map 
\begin{equation}(|A|_s \cup_{|\nabla A|_c} |B|_c)^H \rightarrow X^H\end{equation}
is a weak equivalence and Lemma \ref{keylemma} tells us that 
\begin{equation}|B|_c^H \rightarrow (|A|_s \cup_{|\nabla A|_c} |B|_c)^H\end{equation}
is a weak equivalence.  Therefore the composite is a weak equivalence, and this implies that the adjoint $B \rightarrow S_c(X)$ is a weak equivalence of cyclic sets.

Conversely, assume that the adjoint 
\begin{equation}(A,B, \nabla A \rightarrow B) \rightarrow (S(X^{S^1}),S_c(X), \nabla S(X^{S^1}) \rightarrow S_c(X))\end{equation}
is a weak equivalence.  This implies that $|A|_s \rightarrow X^{S^1}$ is a weak equivalence of simplicial sets and that $|B|_c \rightarrow X$ is a weak equivalence of cyclic sets.  The previous discussion and the ``two out of three'' property for weak equivalences now imply that $|A|_s \cup_{|\nabla A|_c} |B|_c \rightarrow X$ is a weak equivalence.
\end{proof}

\section{Proof of Theorem 1.1}

The proof that $\sC_F\sD$ inherits a model structure from model structures on $\sC$ and $\sD$ when $F$ is Reedy admissible uses the standard technique for lifting model structures to diagram categories indexed by Reedy categories \cite{hovey}, \cite{hirschhorn}.  

\begin{thma}
Let $\sC$ and $\sD$ be model categories and $F : \sC \rightarrow \sD$ be a Reedy admissible functor.  Then $\sC_F\sD$ admits a model structure.  A map $(A,B,FA \rightarrow B) \rightarrow (A^{\prime},B^{\prime},FA^{\prime} \rightarrow B^{\prime})$ is  

\begin{enumerate}
\item{a weak equivalences if $A \rightarrow A^{\prime}$ is a weak equivalence in $\sC$ and $B \rightarrow B^{\prime}$ is a weak equivalence in $\sD$,} 
\item{a fibration if $A \rightarrow A^{\prime}$ is a fibration in $\sC$ and $B \rightarrow B^{\prime}$ is a fibration in $\sD$,}  
\item{a cofibration if $A \rightarrow A^{\prime}$ is a cofibration in $\sC$ and $FA^{\prime} \cup_{FA} B \rightarrow B^{\prime}$ is a cofibration in $\sD$.}
\end{enumerate}
\end{thma}

\begin{proof}
\hspace{5 pt}
\begin{enumerate}
\item{
$\sC_F\sD$ has all small limits and colimits since $F$ preserves colimits and $\sC$ and $\sD$ have all small limits and colimits.
}

\item
{
Weak equivalences satisfy the ``two out of three'' axiom since they do in $\sC$ and $\sD$.
}
\item
{
It is clear that the weak equivalences and fibrations are closed under retracts, since they are defined levelwise.  We need to verify that retracts of cofibrations are cofibrations.  The commutative diagram:

\begin{equation}
\xymatrix{
FA \ar@/^1pc/[rr]\ar[r]\ar[d] & B\ar@/^1pc/[rr]\ar[d] & FC \ar@/^1pc/[rr] \ar[r]\ar[d] & D\ar@/^1pc/[rr]\ar[d] & FA \ar[r]\ar[d] & B\ar[d] \\
FA^{\prime} \ar@/_1pc/[rr] \ar[r] & B^{\prime} \ar@/_1pc/[rr] & FC^{\prime} \ar[r] \ar@/_1pc/[rr] & D^{\prime} \ar@/_1pc/[rr] & FA^{\prime} \ar[r] & B^{\prime}.\\
} 
\end{equation}
\\
implies that $FC^{\prime} \cup_{FC} D \rightarrow D^{\prime}$ is a retract of $FA^{\prime} \cup_{FA} B \rightarrow B^{\prime}$.  Since $FA^{\prime} \cup_{FA} B \rightarrow B^{\prime}$ is a cofibration, we know from the model structure on $\sD$ that $FC^{\prime} \cup_{FC} D \rightarrow D^{\prime}$ is itself a cofibration in $\sD$.  Moreover, it is clear that $C \rightarrow C^{\prime}$ is a cofibration in $\sC$ because it is a retract of $A \rightarrow A^{\prime}$. 
}

\item
{
Now we need to verify the factorization results.  Assume we have a map $(A,B, FA \rightarrow B) \rightarrow (A^{\prime}, B^{\prime}, FA^{\prime} \rightarrow B^{\prime})$.  We will construct a factorization of this map into a trivial cofibration and a fibration (the other case is analogous).  Consider the following diagram:
\begin{equation}
\begin{CD}
FA @>>> B\\
@VVV @VVV \\
FA^{\prime} @>>> B^{\prime}.\\
\end{CD}
\end{equation}
We employ the standard latching space argument.  Choose a factorization of $A \rightarrow A^{\prime}$ as $A \rightarrow C \rightarrow A^{\prime}$ where $A \rightarrow C$ is a trivial cofibration in $\sC$ and $C \rightarrow A^{\prime}$ is a fibration.  This yields a factorization $FA \rightarrow  FC \rightarrow FA^{\prime}$.  So now we have the following diagram:

\begin{equation}
\begin{CD}
FA @>>> B\\
@VVV @VVV \\
FC @>>> ?\\
@VVV @VVV\\
FA^{\prime} @>>> B^{\prime}.\\ 
\end{CD}
\end{equation} 
To complete the diagram choose a factorization of $FC \cup_{FA} B \rightarrow B^{\prime}$ as 
\begin{equation}FC \cup_{FA} B \rightarrow C^{\prime} \rightarrow B^{\prime}\end{equation}
where $C \cup_A B \rightarrow C^{\prime}$ is a trivial cofibration and $C^{\prime} \rightarrow B^{\prime}$ is a fibration, and then put $C^{\prime}$ in for the $?$.  By the assumption on $F$, $B \rightarrow C^{\prime}$ is a weak equivalence.  This yields the factorization 
\begin{equation}(A, B, FA \rightarrow B) \rightarrow (C,C^{\prime}, FC \rightarrow C^{\prime}) \rightarrow (A^{\prime}, B^{\prime}, FA^{\prime} \rightarrow B^{\prime})\end{equation}
in which the first arrow is a trivial cofibration and the second a fibration.
}
\item
{
Finally, we must verify the lifting properties.  Assume we have a trivial cofibration and a fibration (the other case is analogous).  The lifting problem
\begin{equation}
\xymatrix{
(A, A^{\prime}, FA \rightarrow A^{\prime}) \ar[d] \ar[r] & (B, B^{\prime}, FB \rightarrow B^{\prime}) \ar[d] \\
(X, X^{\prime}, FX \rightarrow X^{\prime}) \ar[r] \ar@{.>}[ur] & (Y, Y^{\prime}, FY \rightarrow Y^{\prime})
}
\end{equation}
splits into the following interlocked lifting problems:
\begin{equation}
\xymatrix{
A \ar[r] \ar[d] & B \ar[d] \\
X \ar[r] \ar@{.>}[ur] & Y \\
} 
\qquad
\xymatrix{
A^{\prime} \ar[r] \ar[d] & B^{\prime} \ar[d] \\
X^{\prime} \ar[r] \ar@{.>}[ur] & Y \\
}
\qquad
\xymatrix{
FA \ar@/^1pc/[rr] \ar[r]\ar[d] & A^{\prime}\ar[d]\ar@/^1pc/[rr] & FB \ar[d]\ar[r] & B^{\prime} \ar[d]\\
FX \ar@{.>} [urr] \ar[r]\ar@/_1pc/[rr] & X^{\prime} \ar@/_1pc/[rr] \ar@{.>}[urr] & FY \ar[r] & Y^{\prime}.\\
}
\end{equation}
\\
First, take a lift $X \rightarrow B$ in the lefthand diagram using the model structure on $\sC$.
Now consider the diagram:

\begin{equation}
\begin{CD}
FX \cup_{FA} A^{\prime} @>>> B^{\prime} \\
@VVV @VVV \\
X^{\prime} @>>> Y^{\prime}.\\
\end{CD}
\end{equation}
\\
Here the map $FX \cup_{FA} A^{\prime} \rightarrow B^{\prime}$ is built using the map $FX \rightarrow FB$ obtained from the lift.  Take a lift $X^{\prime} \rightarrow B^{\prime}$ in this diagram using the model structure in $\sD$.  Together, these two lifts provide the desired lifting. 
}
\end{enumerate}
\end{proof}

\begin{remark}
There is a dual version of this result for categories with objects $(A, B, A \rightarrow GB)$ in which $G$ is a co-Reedy admissible functor.  That is, $G$ preserves limits and satisfies an appropriate pullback condition.
\end{remark}

\begin{lemma}
If $\sC$ and $\sD$ are left proper and $F$ is Reedy admissible, then $\sC_F\sD$ is left proper.  If $\sC$ and $\sD$ are right proper and $F$ is a left Quillen functor, then $\sC_F\sD$ is right proper.
\end{lemma}

\begin{proof}
The first assertion follows since fibrations, weak equivalences, and pullbacks are defined levelwise.  For the second assertion, we need that $(A, B, FA \rightarrow B) \rightarrow (A^{\prime}, B^{\prime}, FA^{\prime} \rightarrow B^{\prime})$ a cofibration implies that $B \rightarrow B^{\prime}$ is a cofibration.  If $F$ is a left Quillen functor, this follows from \cite[15.3.11]{hirschhorn}. 
\end{proof}

\section{Acknowledgments}
I wish to express my deep gratitude to Mike Mandell for his generous assistance in solving this problem and to Peter May for his many useful comments and suggestions.  I would also like to thank Phil Hirschhorn and Mark Hovey for helpful discussion and Jan Spalinski for permitting me to reproduce a portion of his thesis in the appendix. 

\appendix

\section{Calculations from Spalinski's thesis}

We reproduce several calculations which appeared in Spalinski's thesis \cite{spalinski:thesis} but not in the paper based on the thesis \cite{spalinski:paper}.  

\subsection{Explicit calculation of $\Psi_r$}

We need to calculate $\Psi_r(\Delta[k])$.  Recall that $\Psi_r$ is the adjoint to $\Phi_r$, where $\Phi_r(X) = (\sd_r X)^{\Z / (r)}$.  It is sufficient to find a cyclic set $A$ such that there is a natural equivalence 
\begin{equation}\hom_{{\bf S}^c}(A,X) \rightarrow \hom_{{\bf S}}(\Delta[k], \Phi_r(X)).\end{equation}
We know that there is an equivalence
\begin{equation}\hom_{{\bf S}}(\Delta[k], \Phi_r(X)) \rightarrow \Phi_r(X)_k\end{equation} given by $f \mapsto f(\iota_k)$.
So it will suffice to exhibit a cyclic set $A$ such that there is a natural equivalence
\begin{equation}\hom_{{\bf S}^c}(A,X) \rightarrow \Phi_r(X)_k.\end{equation}

There is an action of $\Z / (n+1)$ on $|\Lambda[n]|_c$ for $n \geq 1$.  By the Yoneda lemma, each map $\Lambda[n] \rightarrow \Lambda[n]$ is of the form $\hom_{\Lambda^{\op}}(\phi,-)$ for some $\phi : [n] \rightarrow [n] \in \Lambda^{\op}$.  The map corresponding to $t_{n+1}$ has order $n+1$.  This provides the action of $\Z / (n+1)$, and we refer to the generator of this action as $\alpha$.  This action induces an action of $\Z / (n+1)$ on $|\Lambda[n]|_c$. 

\begin{definition}
If $k$ divides $n+1$, let $\Lambda[n\,|\,k]$ denote the orbit space of $\Lambda[n]$ with respect to the action of the subgroup of $\Z / (n+1)$ generated by $\alpha^k$.
\end{definition}

\begin{proposition}
The map
\begin{equation}\hom_{{\bf S}^c}(\Lambda[r(k+1)-1|k+1],X) \rightarrow \Phi_r(X)_k\end{equation}
given by $f \mapsto f[\iota_{r(k+1) - 1}]$ is a bijection and so
$\Psi_r(\Delta[n]) = \Lambda[r(n+1) - 1 | n+1]$.
\end{proposition}

\begin{proof}
First note that the image of $\gamma$ is actually contained in the above fixed point set:
\begin{eqnarray*}
t^{k+1}_{r(k+1)} \cdot \gamma(f) & = & t^{k+1}_{r(k+1)} \cdot f([\iota_{r(k+1)-1}]) \\ & = & f(t^{k+1}_{r(k+1)}[\iota_{r(k+1)-1}]) \\ & = & f([\iota_{r(k+1)-1}]) = \gamma(f).
\end{eqnarray*}
Next, observe that $\gamma$ is onto.  Take $x \in X^H_{r(k+1)-1}$ and consider the map:
\begin{equation}f:\Lambda[r(k+1)-1] \rightarrow X, \quad \iota_{r(k+1)-1} \mapsto x.\end{equation}
Note that $\Im f \subseteq X^H_{r(k+1)-1}$.  Let $z = STD[\iota_{r(k+1)-1}] \in \Lambda[r(k+1)-1]$.  We need to show that $f(\alpha_{k+1}z) = f(z)$.  We have:
\begin{eqnarray*}
f(\alpha_{k+1}z) & = & f(\alpha_{k+1}STD[\iota_{r(k+1)-1}]) = f(STD t^{k+1}_{r(k+1)}[\iota_{r(k+1)-1}]) \\ & = & STD t^{k+1}_{r(k+1)} f[\iota_{r(k+1)-1}] = STD t^{k+1}_{r(k+1)} x = STD x \\ & = & STD f(\iota_{r(k+1)-1}) = f(STD[\iota_{r(k+1)-1}]) = f(z).
\end{eqnarray*}
Hence $f$ factors as:
\begin{equation}\Lambda[r(k+1)-1] \rightarrow \Lambda[r(k+1)-1] / \alpha^{k+1} \rightarrow X.\end{equation}
Here the first map is the quotient map and the second map is $\bar{f}$.  By construction $\gamma(\bar{f}) = x$.  Finally, we need to check that $\gamma$ is injective.  Suppose that $\gamma(f) = \gamma(g)$.  Then $f[\iota_{r(k+1)-1}] = g[\iota_{r(k+1)-1}]$.  Since $[\iota_{r(k+1)-1}]$ generates $\Lambda[r(k+1)-1]/\alpha^{k+1}$, $f = g$.
\end{proof} 

\subsection{Fixed points of $\Psi_r(\Delta[n])$}

The explicit description of $\Psi_r(\Delta[n])$ makes it easy to calculate its $S^1$ fixed points.

\begin{proposition}
For $k \, | \, (n+1)$, $|\Lambda[n \,| \, k]|_c^{S^1} = \emptyset$.
\end{proposition}

\begin{proof}
Let $p \in (\Delta[n] \times S^1)$.  Then we have 
\begin{equation}p = (x_0, x_1, \ldots, x_n, t) \qquad x_i \geq 0, \quad \sum^n_{i=0} x_i = 1, \quad t \in S^1.\end{equation}
Let $\sigma = (0,1,\ldots,n)$, $\tau = \sigma^{-k}$, and $\gamma = e^{2\pi i / \frac{n+1}{k}}$.  The action of $\alpha^k$ on $(\Delta[n] \times S^1)$ is given by 
\begin{equation}\alpha^k(x_0, x_1, \ldots, x_n, t) = (x_{\tau(0)}, x_{\tau(1)}, \ldots, x_{\tau(n)}, \gamma t).\end{equation}
Since $S^1$ is infinite and each orbit of $\alpha^k$ has only finitely many points,
\begin{equation}\{(\Delta[n] \times S^1) / \alpha^k\}^{S^1} = \emptyset.\end{equation}   
\end{proof}

\bibliography{cyclic}

\end{document}